\documentclass[12pt,a4paper]{amsart}
\usepackage[bottom]{footmisc}
\usepackage{mathrsfs}
\usepackage{amsmath}
\usepackage{amsthm}
\usepackage{amsfonts}
\usepackage{latexsym}
\usepackage{graphicx}
\usepackage{hyperref}
\usepackage{amssymb}
\usepackage{subfigure}
\usepackage{epsf}
\usepackage{float}
\usepackage{hyperref}
\usepackage{fancyhdr}
\allowdisplaybreaks \makeatletter
\def\rightharpoonfill@{\arrowfill@\relbar\relbar\rightharpoonup}
\DeclareRobustCommand{\overrightharpoon}{\mathpalette{\underarrow@\rightharpoonfill@}}
\makeatother
\allowdisplaybreaks

\begin{document}
\newcommand{\beq}{\begin{equation}}
\newcommand{\eneq}{\end{equation}}
\newtheorem{thm}{Theorem}[section]
\newtheorem{coro}[thm]{Corollary}
\newtheorem{lem}[thm]{Lemma}
\newtheorem{prop}[thm]{Proposition}
\newtheorem{defi}[thm]{Definition}
\newtheorem{rem}[thm]{Remark}
\newtheorem{cl}[thm]{Claim}
\title{On the extrema of a nonconvex functional with double-well potential in
higher dimensions}
\author{Xiaojun Lu$^{1}$ \ \ \ \ \ \ \ David Yang Gao$^2$}
\pagestyle{fancy}                   
\lhead{X. Lu and D. Y. Gao}
\rhead{Nonconvex functional with double-well potential} 
\thanks{Corresponding author: David Yang Gao}
\thanks{Email addresses:  lvxiaojun1119@hotmail.de(Xiaojun Lu), d.gao@federation.edu.au(David Yang Gao)}
\date{}
\maketitle
\begin{center}
1. Department of Mathematics \& Jiangsu Key Laboratory of
Engineering Mechanics, Southeast University, 210096, Nanjing,
China\\
2. Faculty of Science and Technology, Federation University
Australia, Ballarat, VIC 3350, Australia
\end{center}
\begin{abstract}
This paper mainly addresses the extrema of a nonconvex functional
with double-well potential in higher dimensions through the approach
of nonlinear partial differential equations. Based on the canonical
duality method, the corresponding Euler--Lagrange equation with
Neumann boundary condition can be converted into a cubic dual
algebraic equation, which will help find the local extrema for the
primal problem. In comparison with the 1D case discussed by D. Gao
and R. Ogden, there exists huge difference in higher dimensions,
which will be explained in the theorem.
\end{abstract}
\renewcommand{\abstractname}{R\'{e}sum\'{e}}
\begin{abstract}
Dans cet article, on consid\`{e}re essentiellement les extrema d'un
foncitionnel nonconvexe avec double puits de potentiel dans les
dimensions sup\'{e}rieures. En appliquant la m\'{e}thode de
dualit\'{e} canonique, l'\'{e}quation aux d\'{e}riv\'{e}es
partielles (EDPs) nonlin\'{e}aire, c'est-\`{a}-dire, l'\'{e}quation
d'Euler--Lagrange avec les conditions aux limites de Neumann, peut
\^{e}tre convertie en un dual alg\'{e}brique cubique, qui nous
aidera \`{a} d\'{e}montrer le principe du minimum (ou du maximum)
local pour le probl\`{e}me primal. En comparaison avec le cas 1D
\'{e}tudi\'{e} par D. Gao et R. Ogden, il existe \'{e}norme
diff\'{e}rence dans les dimensions sup\'{e}rieures, et on
l'expliquera dans le th\'{e}or\`{e}me.
\end{abstract}
\section{Introduction}
The double-well potential was first studied by Van der Waals in the
nineteenth century for a compressible fluid whose free energy at a
constant temperature depends on both the density and the density
gradient \cite{waals}. Afterwards, lots of applications of this
nonconvex function have been found in nonlinear sciences, i.e., in
phase transitions of Ericksen's bar \cite{E}, or the mathematical
theory of super-conductivity \cite{G-7}, etc. In this paper, we
consider the fourth-order polynomial defined by\[
H(|\overrightarrow{\gamma}|):={\nu/2}\Big({1/2}|\overrightarrow{\gamma}|^2-\lambda\Big)^2,\
\overrightarrow{\gamma}\in\mathbb{R}^n,\ \nu, \lambda>0\ \text{are
constants,}\
|\overrightarrow{\gamma}|^2=\overrightarrow{\gamma}\cdot\overrightarrow{\gamma}.\]
In quantum mechanics, if $\overrightarrow{\gamma}$ represents Higgs'
field strength, then $H(|\overrightarrow{\gamma}|)$ is the Higgs'
potential \cite{higgs}. It was discovered in the context of
post-buckling analysis \cite{gao-mrc96} that the stored potential
energy of a large deformed beam model in 1D is exactly a double-well
function, where each potential well represents a possible buckled
beam state, and the local maximizer is corresponding to the
unbuckled state \cite{G-10}. As a matter of fact, the polynomial is
also the well-known Landau's second-order free energy, each of its
local minimizers represents a possible phase state of the material,
while each local maximizer characterizes the critical conditions
that lead to the phase transitions etc. \cite{G-6,G-7}.\\

The purpose of this paper is to find the extrema of the following
nonconvex total potential energy functional in higher dimensions,
\begin{equation}
I[u]:=\int_\Omega \Big(H(|\nabla u|)-fu\Big)dx, \label{eq-pp}
\end{equation}
where $\Omega={\rm Int}\Big\{\mathbb{B}(O,R_1)\backslash
\mathbb{B}(O,R_2)\Big\}$, $R_1>R_2>0$, $\mathbb{B}(O,R_1)$ and
$\mathbb{B}(O,R_2)$ denote two open balls with center $O$ and radii
$R_1$ and $R_2$ in the Euclidean space $\mathbb{R}^n$, respectively.
``Int" denotes the interior points. In addition, let
$\Sigma_1:=\{x:|x|=R_1\}$, and $\Sigma_2:=\{x:|x|=R_2\}$, then the
boundary $\partial\Omega=\Sigma_1\cup\Sigma_2$. The radially
symmetric function $f\in C(\overline{\Omega})$ satisfies the
normalized balance condition \beq\int_\Omega f(|x|)dx=0,\eneq and
\beq f(|x|)=0\ \text{if and only if}\ |x|=R_3\in(R_2,R_1).\eneq
Moreover, its $L^1$-norm is sufficiently small such that \beq
\|f\|_{L^1(\Omega)}< 4\lambda\nu
R_2^{n-1}\sqrt{2\lambda\pi^n}/(3\sqrt{3}\Gamma(n/2)),\eneq where
$\Gamma$ stands for the Gamma function. This assumption is
reasonable since large $\|f\|_{L^1(\Omega)}$ may possibly lead to
instant fracture. The deformation $u$ is subject to the following
three constraints, \beq u\ \text{is radially symmetric on}\
\overline{\Omega},\eneq \beq u\in W^{1,\infty}(\Omega)\cap
C(\overline{\Omega}),\eneq \beq \nabla u\cdot\overrightarrow{n}=0\
\text{on both}\ \Sigma_1\ \text{and}\ \Sigma_2, \eneq where
$\overrightarrow{n}$ denotes the unit outward normal on
$\partial\Omega$.

By variational calculus, one derives a correspondingly nonlinear
Euler--Lagrange equation for the primal nonconvex functional,
namely, \beq \displaystyle {\rm div}\Big(\nabla H(|\nabla
u|)\Big)+f=0 \ \text{\rm in}\ \Omega, \eneq equipped with the
Neumann boundary condition (7). Clearly, (8) is a highly nonlinear
partial differential equation which is difficult to solve by the
direct approach or numerical method \cite{JH,LIONS}. However, by the
canonical duality method, one is able to demonstrate the existence
of solutions for this type of equations.\\

This paper is aimed to solve the challenging nonconvex variational
problem by using the {\it canonical duality theory}, which can be
applied in solving a large class of nonconvex/nonsmooth/discrete
problems in multidisciplinary fields including mathematical physics,
global optimization, computational science, etc. For instance, Gao
and Ogden have introduced this method to the 1D problems in finite
deformation mechanics \cite{G-5} and phase transitions of the
Ericksen's bar \cite{G-6}. Their work showed that the nonlinear
differential equations can be converted into dual algebraic
equations and help us find the possible nonsmooth solutions.\\

Before introducing the main result, we denote
\[ F(r):=-1/r^n\int^r_{R_2}f(\rho)\rho^{n-1} d\rho.\
r\in[R_2,R_1].
\]
Next, we define a polynomial of third order as follows, \[
E(y):=\displaystyle 2y^2(\lambda+y/\nu), y\in[-\nu\lambda,+\infty).
\]
Furthermore, for any $A\in[0,8\lambda^3\nu^2/27)$,
\[E_3^{-1}(A)\leq E_2^{-1}(A)\leq E_1^{-1}(A)\]
stand for the three real-valued roots for the equation $E(y)=A.$\\

At the moment, we would like to introduce the theorem of multiple
extrema for the nonconvex functional (2).
\begin{thm}
For any radially symmetric function $f\in C(\overline{\Omega})$
satisfying (2)--(4), we have three solutions for the nonlinear
Euler--Lagrange equation (8) equipped with the Neumann boundary
condition, namely
\begin{itemize}
\item
For any $r\in[R_2,R_1]$, $\bar{u}_1$ defined below is a local
minimizer for the nonconvex functional (2). \beq
\bar{u}_1(|x|)=\bar{u}_1(r):=
\int^{r}_{R_2}F(\rho)\rho/E_1^{-1}(F^2(\rho)\rho^2)d\rho+C_1,\
\forall\ C_1\in\mathbb{R}.\eneq
\item For any $r\in[R_2,R_1]$, $\bar{u}_2$ defined below is a local
minimizer for the nonconvex functional (2) in 1D. While for the
higher dimensions $n\geq 2$, $\bar{u}_2$ is not necessarily a local
minimizer for (2) in comparison with the 1D case. \beq
\bar{u}_2(|x|)=\bar{u}_2(r):=
\int^{r}_{R_2}F(\rho)\rho/E_2^{-1}(F^2(\rho)\rho^2)d\rho+C_2,\
\forall\ C_2\in\mathbb{R}.\eneq
\item For any $r\in[R_2,R_1]$, $\bar{u}_3$ defined below is a local
maximizer for the nonconvex functional (2). \beq
\bar{u}_3(|x|)=\bar{u}_3(r):=
\int^{r}_{R_2}F(\rho)\rho/E_3^{-1}(F^2(\rho)\rho^2)d\rho+C_3,\
\forall\ C_3\in\mathbb{R}. \eneq
\end{itemize}
\end{thm}
\begin{rem}
From the above theorem, one knows, it is incorrect to simply
generalize the 1D case discussed in \cite{G-6} to higher dimensions.
Each of these solutions is a critical point of $I$, i.e., it could
be either an extremum or a saddle point of the total potential. This
phenomenon has been verified by Ericksen who proved that many local
solutions are metastable and may have arbitrary number of phase
interfaces.
\end{rem}
\begin{rem} Compared with
convex problems, a fundamentally different issue in nonconvex
analysis is that the solutions of the boundary-value problem is not
equivalent to the associated minimum variational problem and it is
extremely difficult to use traditional direct approaches in solving
the nonconvex variational problems. In particular, it is discovered
that for certain given external loads, a global or local minimizer
is nonsmooth and cannot be determined by any Newton-type numerical
methods.
\end{rem}
The rest of the paper is organized as follows. In Section 2, first
we introduce some useful notations which will simplify our proof
considerably. Then, we apply the canonical dual transformation to
deduce a perfect dual problem corresponding to the primal nonconvex
variational problem and a pure complementary energy principle. In
the final analysis, we apply the canonical duality theory to prove
Theorem 1.1.
\section{Proof of the main result}
\subsection{Some useful notations}

\begin{itemize}
\item $\overrightarrow{\sigma}$ is the G\^{a}teaux derivative of
$H$ given by
\[
\overrightarrow{\sigma}(x)=(\sigma^{(1)}(x),\cdots,\sigma^{(n)}(x))=\nu(1/2|\nabla
u|^2-\lambda)\nabla u.
\]
From the physician' viewpoint, if $\nabla u$ represents a
deformation gradient, then the vector $\overrightarrow{\sigma}$ is
the so-called first Piola-Kirchhoff stress in the finite deformation
theory.
\item $\Phi$ is a nonlinear
geometric mapping defined as
\[
\Phi(u):=1/2|\nabla u|^2.
\]
For convenience's sake, denote $\xi:=\Phi(u).$ It is evident that
$\xi$ belongs to the function space $\mathscr{U}$ given by
\[
\mathscr{U}:= \Big\{\phi\in L^\infty(\Omega)\Big| \phi\geq0\Big\}.
\]
\item $\Psi$ is a canonical energy
defined as
\[
\Psi(\xi):=\nu/2(\xi-\lambda)^2,
\]
which is a convex function with respect to $\xi$. For simplicity,
denote $\zeta:=\nu(\xi-\lambda)$, which is the G\^{a}teaux
derivative of $\Psi$ with respect to $\xi$. Moreover, $\zeta$ is
invertible with respect to $\xi$ and belongs to the function space
$\mathscr{V}$,
\[\mathscr{V}:=\Big\{\phi\in
L^\infty(U)\Big| \phi\geq-\nu\lambda\Big\}.
\]
\item
$\Psi_\ast$ is defined as
\[
\Psi_\ast(\zeta):=\xi\zeta-\Psi(\xi)=\zeta^2/(2\nu)+\lambda\zeta.
\]
\end{itemize}
\subsection{Canonical duality techniques}
\begin{defi}
By Legendre transformation, one defines a Gao--Strang total
complementary energy functional $\Xi$,
\[
\Xi(u,\zeta):=\displaystyle\int_{\Omega}\Big\{\Phi(u)\zeta-\Psi_\ast(\zeta)
-fu\Big\}dx.
\]
\end{defi}
Next we introduce an important {\it criticality criterium} for the
Gao-Strang total complementary energy functional.
\begin{defi}
$(\bar{u}, \bar{\zeta})$ is called a critical pair of $\Xi$ if and
only if \beq D_{u}\Xi(\bar{u},\bar{\zeta})=0, \eneq \beq
D_{\zeta}\Xi(\bar{u},\bar{\zeta})=0, \eneq where $D_{u}, D_{\zeta}$
denote the partial G\^ateaux derivatives of $\Xi$, respectively.
\end{defi} Indeed, by variational calculus, one has the following
observation from (12) and (13).
\begin{lem}
On the one hand, for any fixed $\zeta\in\mathscr{V}$, $(12)$ is
equivalent to the equilibrium equation
\[
\begin{array}{ll}\displaystyle {\rm div}(\zeta \nabla\bar{u})+f=0& \
\text{\rm in}\ \Omega,\end{array}
\]
with the Neumann boundary condition. On the other hand, for any
fixed $u$ subject to (5)-(7), (13) is consistent with the notations
before,
\[
\Phi(u)=D_{\zeta}\Psi_\ast(\bar{\zeta}).
\]
\end{lem}
Lemma 2.3 indicates that $\bar{u}$ from the critical pair
$(\bar{u},\bar{\zeta})$ solves the Euler--Lagrange equation (8).
\begin{defi}
From Definition 2.1, one defines the Gao--Strang pure complementary
energy $I_d$ in the form
\[
I_d[\zeta]:=\Xi(\bar{u},\zeta),
\]
where $\bar{u}$ solves the Euler--Lagrange equation (8).
\end{defi}
To simplify the discussion, one uses another representation of the
pure energy $I_d$ given by the following lemma through integrating
by parts.
\begin{lem} The
pure complementary energy functional $I_d$ can be rewritten as \beq
I_d[\zeta]=-1/2\int_{\Omega}\Big\{{|\overrightarrow{\sigma}|^2/\zeta}+2\lambda\zeta+\zeta^2/\nu\Big\}dx,
\eneq where $\overrightarrow{\sigma}$ satisfies \beq {\rm
div}\overrightarrow{\sigma}+f=0\ \text{in}\ \Omega, \eneq equipped
with $\overrightarrow{\sigma}\cdot\overrightarrow{n}=0$ on
$\partial\Omega$.
\end{lem}
With the above discussion, indeed, by calculating the G\^{a}teaux
derivative of $I_d$ with respect to $\zeta$, one has \begin{lem} The
variation of $I_d$ with respect to $\zeta$ leads to the cubic dual
algebraic equation (DAE), namely \beq
|\overrightarrow{\sigma}|^2=2\bar{\zeta}^2(\lambda+\bar{\zeta}/\nu),
\eneq where $\bar{\zeta}$ is from the critical pair
$(\bar{u},\bar{\zeta})$.
\end{lem}
\begin{rem}
From (16), it is easy to check that $|\overrightarrow{\sigma}|^2$
has a maximum $8\lambda^3\nu^2/27$ at $\bar{\zeta}=-2\lambda\nu/3$
and a minimum 0 at $\bar{\zeta}=0$.
\begin{itemize}
\item If $|\overrightarrow{\sigma}|^2=8\lambda^3\nu^2/27$, then there exist
two real roots.
\item If
$|\overrightarrow{\sigma}|^2\in(8\lambda^3\nu^2/27,\infty)$, then
there exists only one positive real root.
\item If
$|\overrightarrow{\sigma}|^2\in[0,8\lambda^3\nu^2/27)$, then there
exist three real roots listed below, \beq
\bar{\zeta}_1>0>\bar{\zeta}_2>-2\nu\lambda/3>\bar{\zeta}_3>-\nu\lambda.
\eneq
\item For $|\overrightarrow{\sigma}|^2=0$, there also exist three real
roots such as \beq \bar{\zeta}_1=\bar{\zeta}_2=0,\ \
\bar{\zeta}_3=-\nu\lambda. \eneq
\end{itemize}
\end{rem}
\subsection{Proof of Theorem 1.1}
Actually, a radially symmetric solution for the Euler--Lagrange
equation (8) is of the form
\[
\overrightarrow{\sigma}=F(r)(x_1,\cdots,x_n)=F\Big(\sqrt{\sum_{i=1}^nx_i^2}\Big)(x_1,\cdots,x_n),
\]
where $F$ is the unique solution for the nonhomogeneous linear
differential equation
\[
F'(r)+nF(r)/r=-f(r)/r,\ \ \ r\in[R_2,R_1].
\]
with $F(R_2)=0$. Furthermore, the normalized balance condition (2)
assures that $F(R_1)=0$, which indicates
$\overrightarrow{\sigma}\cdot\overrightarrow{n}=0$ on
$\partial\Omega$.\\

From the above discussion, one deduces that once
$\overrightarrow{\sigma}$ is given, then an analytic solution of the
Euler--Lagrange equation (8) can be presented as \beq
\bar{u}_i(|x|)=\displaystyle\int^{x}_{x_0}\overrightarrow{\eta_i}\overrightarrow{dt},
\eneq where $x\in \overline{\Omega}, x_0\in\partial \Omega$,
$\overrightarrow{\eta_i}=(\eta_i^{(1)},\eta_i^{(2)},\cdots,\eta_i^{(n)}):=\overrightarrow{\sigma}/\bar{\zeta}_i$,
which satisfies the condition for path independent integrals,
namely, for $i=1,2,3$,
\[
\partial_{x_j}\eta_{i}^{(k)}-\partial_{x_k}\eta_{i}^{(j)}=0,\
\ j,\ k=1,\cdots,n.
\]
By direct calculation for $I$ and $I_d$, respectively, it is easy to
check that the pure complementary energy functional $I_d$ is
perfectly dual to the total potential energy functional $I$, and the
identities
\[I[\bar{u}_i]=I_d[\bar{\zeta}_i],\ \ i=1,2,3\] indicate there is no duality gap between the primal and dual
variational problems.\\

Next, we prove the local extrema. On the one hand, for any test
function $\phi\in W^{1,\infty}(\Omega)$, the second variational form
$\delta_\phi^2I$ is equal to\beq \nu\int_\Omega\Big\{|\nabla
\bar{u}\cdot\nabla\phi|^2+\Big(1/2|\nabla
\bar{u}|^2-\lambda\Big)|\nabla\phi|^2\Big\}dx.\eneq On the other
hand, for any test function $\psi\in\mathscr{V}$, the second
variational form $\delta_\psi^2I_d$ is equal to
\beq-\int_\Omega\Big\{{|\overrightarrow{\sigma}|^2}/{\bar{\zeta}^3}+{1}/{\nu}\Big\}\psi^2dx.
\eneq  Actually, (3) and (4) indicate $0<
F^2(r)r^2<8\lambda^3\nu^2/27$ for any $r\in(R_2,R_1)$. Indeed, let
\[G(r):=\int_{R_2}^rf(\rho)\rho^{n-1}d\rho.\]
Its derivative is $G'(r)=f(r)r^{n-1}$ and $G(R_2)=G(R_1)=0$. From
(3), since $f(|x|)=0$ if and only if $|x|=R_3\in(R_2,R_1)$, without
loss of generality, we assume that
\[
f(r)=\left\{\begin{array}{lll}
>0,& r\in[R_2,R_3);\\
\\
<0, & r\in(R_3,R_1].
\end{array}\right.\]
From the above assumption and the fact $G(R_2)=G(R_1)=0$, one has
$G(r)>0$, $r\in(R_2,R_1)$. As a result, $F(r)$ does not change its
sign in $(R_2,R_1)$ and
\[F^2(r)r^2=r^{2-2n}G^2(r)>0,\ r\in(R_2,R_1).\] From (4), one has
\[
\begin{array}{lll}F^2(r)r^2&=&r^{2-2n}\Big(\int_{R_2}^rf(\rho)\rho^{n-1}d\rho\Big)^2\\
\\
&\leq&r^{2-2n}\Big(\int_{R_2}^r|f(\rho)|\rho^{n-1}d\rho\Big)^2\\
\\
&=&r^{2-2n}\Gamma^2(n/2)/(4\pi^{n})\|f\|^2_{L^1(\Omega)}\\
\\
&<&(R_2/r)^{2n-2}8\lambda^3\nu^2/27\\
\\
&<& 8\lambda^3\nu^2/27.
\end{array}
\]

\begin{itemize}
\item According to the definition of $\zeta$, one knows immediately
that for $\bar{\zeta}_1>0$,
\[\delta^2_\phi I(\bar{u}_1)\geq 0,\ \
\delta^2_\psi I_d(\bar{\zeta}_1)\leq0.\]
\item Since $\bar{\zeta}_3<-2/3\nu\lambda$, then
$\delta_\phi^2I_d(\bar{\zeta}_3)\leq0$ and $\lambda>3/2|\nabla
\bar{u}_3|^2$. In this case,
\begin{eqnarray*}
\delta^2_\phi I(\bar{u}_3)& \leq & \nu\int_\Omega\Big\{|\nabla
\bar{u}_3|^2|\nabla\phi|^2+\Big(1/2|\nabla
\bar{u}_3|^2-\lambda\Big)|\nabla\phi|^2\Big\}dx \\
& = & \nu\int_\Omega\Big(3/2|\nabla
\bar{u}_3|^2-\lambda\Big)|\nabla\phi|^2dx\leq0.
\end{eqnarray*}
\item Indeed, $-2\nu\lambda/3<\bar{\zeta}_2<0$ indicates
$\lambda\in(1/2|\nabla \bar{u}_2|^2,3/2|\nabla \bar{u}_2|^2)$. In
this case, $\delta^2_\psi I_d(\bar{\zeta}_2)\geq0$. As for the 1D
case, $\delta^2_\phi I(\bar{u}_2)\geq 0$, which induces that
$\bar{u}_2$ is a local minimizer in 1D. But for the higher
dimensional case $n\geq2$, $\bar{u}_2$ is not necessarily a local
minimizer, which is in direct contrast to the 1D case.
\end{itemize}
In the final analysis, together with the fact $F\in C^1[R_2,R_1]$,
definition of $E$, (17) and (18), we have
\[\lim_{\rho\to R_2^+}F(\rho)\rho/E_i^{-1}(F^2(\rho)\rho^2)<\infty,\ i=1,2,3,\]
and
\[\lim_{\rho\to R_1^-}F(\rho)\rho/E_i^{-1}(F^2(\rho)\rho^2)<\infty,\ i=1,2,3,\]
which indicate $\bar{u}_i\in C(\overline{\Omega}),\ i=1,2,3$. And
the proof of Theorem 1.1 is
concluded.\\

{\bf Acknowledgment}: The main results in this paper were obtained
during a research collaboration at the Federation University
Australia in August, 2016. The first author wishes to thank
Professor David Y. Gao for his hospitality and financial support.
This project is partially supported by the US Air Force Office of
Scientific Research (AFOSR FA9550-10-1-0487), Natural Science
Foundation of Jiangsu Province (BK 20130598), National Natural
Science Foundation of China (NSFC 71273048, 71473036, 11471072), the
Scientific Research Foundation for the Returned Overseas Chinese
Scholars. This work is also supported by Open Research Fund Program
of Jiangsu Key Laboratory of Engineering Mechanics, Southeast
University (LEM16B06). In particular, the authors also express their
deep gratitude to the referees for their careful reading and useful
remarks.

\end{document}